\documentclass[12pt]{article}

\usepackage{amstext    }
\usepackage{amsthm    }
\usepackage{a4}
\usepackage[mathscr]{eucal} 
\usepackage{mathrsfs}

\usepackage{amsmath}
\usepackage{amssymb}
\usepackage{amscd}

\numberwithin{equation}{section}

\newtheorem{theorem}{Theorem}[section]
\newtheorem{definition}[theorem]{Definition}
\newtheorem{proposition}[theorem]{Proposition}
\newtheorem{corollary}[theorem]{Corollary}
\newtheorem{lemma}[theorem]{Lemma}
\newtheorem{remark}[theorem]{Remark}

\newtheorem{example}[theorem]{Example}

\newcommand{\cali}[1]{\mathscr{#1}}

\newcommand{\volume}{{\rm volume}}
\newcommand{\lov}{{\rm lov}}
\newcommand{\supp}{{\rm supp}}

\newcommand{\dist}{{\rm distance}}

\newcommand{\id}{{\rm id}}

\newcommand{\Bc}{\cali{B}}

\newcommand{\Fc}{\cali{F}}

\newcommand{\Lc}{\cali{L}}

\newcommand{\Oc}{\cali{O}}

\newcommand{\FS}{{\rm FS}}

\newcommand{\C}{\mathbb{C}}
\newcommand{\N}{\mathbb{N}}

\renewcommand{\P}{\mathbb{P}}

\title{Upper bound for topological entropy of a meromorphic correspondence}
\author{Tien-Cuong Dinh and Nessim Sibony}

\begin{document}
\maketitle

\begin{abstract}
Let $f$ be a meromorphic correspondence on a compact K{\"a}hler manifold.
We show that the topological entropy of $f$ is bounded from above by
the logarithm of its maximal dynamical degree. An analogous estimate
for the entropy on subvarieties is given. We also discuss a notion
of Julia and Fatou sets.
\end{abstract}

\noindent
{\bf AMS classification :} 37A35, 37B40, 37F05, 32U40

\noindent
{\bf Key-words :} correspondence, entropy,
  dynamical degree, Julia and Fatou sets

\section{Introduction} \label{section_introduction}

Let $(X,\omega)$ be a compact K{\"a}hler manifold of dimension $k$. A
meromorphic correspondence $f:X\rightarrow X$ is a meromorphic multivalued self-map
on $X$. The precise definition will be given in Section
\ref{section_correspondence}. 
One can compose correspondences and consider the dynamical system
associated to $f$, i.e. study the sequence of iterates 
$f^n:=f\circ\cdots \circ f$, $n$ times, of $f$. 
Any projective manifold admits dynamically interesting correspondences.
The topological entropy $h(f)$ of
$f$ is defined as in \cite{Bowen, Gromov, Friedland}, see 
Section \ref{section_entropy}. It measures the divergence of the orbits of $f$
and the complexity of the associated dynamical system.

In this paper, we show that $h(f)$ is
bounded from above by the logarithm of the maximal dynamical degree
of $f$ which is easier to compute or estimate. 
The dynamical degree $d_p(f)$ of order $p$ measures the growth of
the norms of $f^n$ acting on the cohomology group $H^{p,p}(X,\C)$ when
$n$ tends to infinity, see
Section \ref{section_correspondence}. 
Let $\Gamma_{[n]}$ denote the graph of $(f,f^2,\ldots,f^n)$ in
$X^{n+1}$.
We will use the following intermediate indicator,
introduced by Gromov \cite{Gromov}, which measures the growth of the
volume of $\Gamma_{[n]}$ :
$$\lov(f):=\lim_{n\rightarrow \infty} \log \left(\volume(\Gamma_{[n]})^{1/n}\right).$$
We will see that the last limit always exists.
Our main result is the following theorem which is new even for
holomorphic correspondences.
It answers a problem raised  by Gromov \cite{Gromov, Friedland}. 

\begin{theorem} \label{th_entropy}
Let $f$ be a meromorphic correspondence on a compact K{\"a}hler manifold $(X,\omega)$ of dimension $k$.
Let $d_p(f)$ denote the dynamical degree of order $p$ of $f$.
Then 
$$h(f)\leq \lov(f)=\max_{0\leq p\leq k}\log d_p(f).$$
\end{theorem}
The case of holomorphic maps was proved by Gromov \cite{Gromov}, see
also \cite{Friedland}, and  the case of meromorphic maps was proved by the authors
in \cite{DinhSibony2, DinhSibony3}. For other contexts, see \cite{DinhSibony1,
  DinhSibony4, deThelin, Dinh2}
and the references therein. 
The proofs in the previous cases cannot be extended to
correspondences. We need here new geometric ingredients. 
In the last two sections we extend the previous result to the entropy
of $f$ on a subvariety of $X$ and we discuss a notion of Julia and
Fatou sets for correspondences. Our goal is also to develop a calculus
for meromorphic correspondences.

Note that if $f$ is a holomorphic self-map on $X$, by Yomdin's 
theorem \cite{Yomdin}, we have $h(f)\geq\max_p\log d_p(f)$; then
 $h(f)=\max_p\log d_p(f)$, see also \cite{MisiurewiczPrzytycki,
   Przytycki, Newhouse, FornaessSibony, Sibony}. 
However, this is false for holomorphic
correspondences, even in dimension 1.
Let $(z,w)$ denote the canonical affine coordinates of $\C\times\C$ in $\P^1\times
\P^1$. Consider the correspondence $f$ on $\P^1$
with irreducible graph $\Gamma\subset \P^1\times\P^1$ of equation
$w^2=z^2+1$. The reader can check easily that 
$d_0(f)=d_1(f)=2$ and $h(f)=0$.


\section{Regularization of currents} \label{section_current} 

Recall that the mass of a positive $(p,p)$-current $T$ on a compact
K{\"a}hler manifold $(X,\omega)$ of dimension $k$ 
is given by $\|T\|:=\langle T,\omega^{k-p}\rangle$. It depends
continuously on $T$.
When $T$ is positive closed, $\|T\|$ depends only on the class of $T$
in $H^{p,p}(X,\C)$.
In order to simplify the
notation, if $Y$ is an analytic set of pure
dimension in $X$, we often denote by $Y$, instead of $[Y]$, the current of integration on $Y$ and
by $\|Y\|$ its mass. 
The main tool used in the proof of Theorem \ref{th_entropy} is
the following result.

\begin{theorem}[\cite{DinhSibony2, DinhSibony3}] \label{th_regularization}
Let $(X,\omega)$ be a compact K{\"a}hler manifold of dimension $k$. Let
$T$ be a positive closed $(p,p)$-current on $X$. Then there are
positive closed $(p,p)$-currents $T^\pm$ and a constant $c>0$
independent of $T$ such that
\begin{enumerate}
\item[i)] $T=T^+-T^-$ and $\|T^\pm\|\leq c\|T\|$;
\item[ii)] $T^\pm$ are limits of smooth positive closed
  $(p,p)$-forms on $X$. 
\end{enumerate}
\end{theorem}

We deduce the following consequence.

\begin{corollary} \label{cor_intersection}
Let $\pi:(X_1,\omega_1)\rightarrow (X_2,\omega_2)$ be a holomorphic
map between two compact K{\"a}hler
manifolds.
Let $Y\subset X_1$ be an analytic subset of pure dimension 
and let $Y'$ be a Zariski open subset of $Y$ such that the
restriction $\tau$ of $\pi$ to $Y'$ is locally a submersion on $X_2$. If 
$T$ is a positive closed current on $X_2$, then $\tau^*(T)$ extends to
 a positive closed current on $X_1$ such that
$$\|\tau^*(T)\|\leq c\|Y\|\ \|T\|,$$
where the constant $c>0$ depends only on $(X_1,\omega_1)$,
$(X_2,\omega_2)$ and $\pi$.
\end{corollary}
\proof
Observe that $\tau^*(T)$ defines a positive closed current on the
Zariski open subset $X\setminus(Y\setminus Y')$ of $X$. If $\tau^*(T)$ has finite mass, a theorem of
Skoda \cite{Skoda} implies that its trivial extension defines a
positive closed current on $X$. Then, we only need to estimate the
mass of $\tau^*(T)$.

The constants that we use here are independent of $Y$, $Y'$ and $T$.
By Theorem \ref{th_regularization}, there are smooth positive closed
forms $T_n$ converging to a current $T'\geq T$ and such that
$\|T_n\|\leq c' \|T\|$.
It follows that there is a constant $c''>0$ such that
$c''\|T\|\omega^p_2-T_n$ is cohomologous to a smooth positive closed form for
every $n$. Here, $(p,p)$ is the bidegree of $T$ and we
use the fact that $H^{p,p}(X_2,\C)$ has finite dimension. 
So the class of $T_n$ is bounded by the class of
$c''\|T\|\omega_2^p$. Since $\tau$ is locally a
submersion on $Y'\subset Y$, we have 
$$\|\tau^*(T)\|\leq \limsup_{n\rightarrow\infty}
\|Y\wedge \pi^*(T_n)\|\leq c''\|T\|\ \|Y\wedge \pi^*(\omega_2^p)\|\leq
c\|Y\| \ \|T\|.$$
In the last inequalities, we use the fact that the mass of a positive closed current
depends only on its cohomology class.
\endproof

\begin{remark} \label{rk_intersection} \rm
If $T$ is the current of integration on a subvariety $Y_2\subset X_2$
then we obtain from the previous corollary that
$$\| Y'\cap \pi^{-1}(Y_2)\|\leq c\|Y\|\ \|Y_2\|.$$
This is a B{\'e}zout type theorem in which
we do not assume that the intersection $Y\cap\pi^{-1}(Y_2)$ is of pure dimension.
\end{remark}


\section{Correspondences and dynamical degrees} \label{section_correspondence}

Let $\pi_1$ and $\pi_2$ denote the canonical projections of $X^2$ 
onto its factors.
A {\it meromorphic correspondence} $f$ on $X$ 
is given by a finite holomorphic chain 
$\Gamma=\sum\Gamma_i$
such that
\begin{enumerate}
\item[{\it i})] for each $i$, $\Gamma_i$ is an irreducible analytic 
subset of dimension $k$ of $X^2$;
\item[{\it ii})] $\pi_1$ and $\pi_2$ restricted to each $\Gamma_i$ are surjective.   
\end{enumerate}
We call $\Gamma$ {\it the graph} of $f$.
We do not assume that the $\Gamma_i$'s are smooth or distinct. Of course, we can write $\Gamma=\sum n_j\Gamma_j'$
where $n_j$ are positive integers and $\Gamma_j'$ are distinct irreducible analytic sets.
Then, a generic point in the support $\cup\Gamma_j'$ of $\Gamma$ belongs to a unique $\Gamma_j'$ and $n_j$ is 
called  {\it the multiplicity} of $\Gamma$ at $x$. 
In what follows we use the notation $\sum\Gamma_i$. The indice $i$
permits to count the multiplicities. Let $\Gamma^{-1}$ denote the
symmetric of $\Gamma$ with respect to the diagonal of $X^2$. The
correspondence $f^{-1}$ associated to $\Gamma^{-1}$ is called
{\it the adjoint} of $f$. 

Observe that if $\Omega$ and $\Omega'$ are dense Zariski open sets in
$X$, then, by condition {\it ii}), all components of $\Gamma$ intersect
$\pi_1^{-1}(\Omega)\cap\pi_2^{-1}(\Omega')$. Hence, $\Gamma$ is the
closure of its restriction to
$\pi_1^{-1}(\Omega)\cap\pi_2^{-1}(\Omega')$.
We will use this property several times.

We define formally $f:=\pi_2\circ(\pi_{1|\Gamma})^{-1}$. More precisely,
if $A$ is a subset of $X$, define
$$f(A):=\pi_2(\pi_1^{-1}(A)\cap\Gamma) \ \mbox{ and }\ f^{-1}(A)=
\pi_1(\pi_2^{-1}(A)\cap\Gamma).$$
So, generically the fibers $f(x)$ and $f^{-1}(x)$ are finite subsets of $X$.
The sets 
$$I_1(f):=\big\{x\in X,\ \dim \pi_1^{-1}(x)\cap\Gamma >0\big\}$$
and
$$I_2(f):=\big\{x\in X,\ \dim \pi_2^{-1}(x)\cap\Gamma >0\big\}$$
are {\it the first and second indeterminacy sets} of $f$; they are of
codimension $\geq 2$.
One can compare the restriction of $\pi_1$ to $\Gamma$ with a blow up
of $X$ along $I_1(f)$ and $\pi_1^{-1}(I_1(f))\cap\Gamma$ is contracted
by $\pi_1$ to
$I_1(f)$. 
If $I_1(f)=\varnothing$ we say that $f$ 
is {\it holomorphic}.
If generic fibers of $\pi_{1|\Gamma}$ contain only one point, we 
obtain a dominant meromorphic self-map on $X$.

We can compose correspondences.  Let $f$ and $f'$ be two correspondences
on $X$ of graphs $\Gamma=\sum_i\Gamma_i$ and
$\Gamma'=\sum_j\Gamma_j'$ in $X^2$. 
Then, the graph of $f\circ f'$ is 
equal to $\Gamma\circ\Gamma':=\sum_{i,j}\Gamma_i\circ\Gamma_j'$, where $\Gamma_i\circ \Gamma_j'$ 
is defined as follows. 

Let $P_i(f)$ denote the smallest analytic subset of $X$ such that $\pi_i$ restricted 
to $\Gamma\setminus \pi_i^{-1}(P_i(f))$ defines an unramified  covering over $X\setminus P_i(f)$. 
Let $\Omega\subset X\setminus P_1(f)$ be a dense Zariski open subset of $X$.
Let $\Omega'\subset X\setminus P_1(f')$ be a similar Zariski open set for $f'$
such that $f'(\Omega')\subset\Omega$. We can choose for example
$\Omega'=(X'\setminus P_1(f'))\setminus f'^{-1}(X\setminus\Omega)$.
Let $\Sigma$ be the closure in $X^2$ of the set  
$$\big\{(x,z)\in \Omega'\times X,\quad \mbox{there is } 
y\in X \mbox{ with } (x,y)\in\Gamma_j' \mbox{ and }
(y,z)\in\Gamma_i\big\}.
\footnote{if we take $(x,z)\in X\times X$, we may obtain components
  whose projections are not equal to $X$; this is the case for the
  iterates of
  non algebraically stable maps in the sense of \cite{Sibony}}  $$
The composition $\Gamma_i\circ\Gamma_j'$ is the holomorphic $k$-chain with support in $\Sigma$ where 
the multiplicity of a generic point $(x,z)$ is defined as  
the number of $y$'s satisfying the previous conditions; quite
generically the multiplicity is one. 

Observe that $\Sigma$ and $\Gamma\circ\Gamma'$ do not depend on the choice of $\Omega$ and $\Omega'$. 
Note that
compositions of irreducible correspondences can be reducible.  This is the reason why we have to deal with 
multiplicities. 
For example, 
if the graph $\Gamma$ of an irreducible correspondence $f$ is symmetric with respect to the diagonal of $X^2$ and
if the degree of $\pi_{i|\Gamma}$ is larger than 1 then $f^2$ is
reducible since its graph contains the diagonal $\Delta$ 
of $X^2$ as one component. Note also that the graph of $f\circ f^{-1}$
contains $\Delta$ but in general we do not have
$f\circ f^{-1}=\id$.  

Correspondences act on smooth forms. If $\alpha$ is a smooth
$(p,p)$-form on $X$, define
$$f^*(\alpha):=(\pi_1)_*\big(\Gamma\wedge \pi_2^*\alpha\big)\quad
\mbox{and}\quad f_*(\alpha):=(\pi_2)_*\big(\Gamma\wedge
\pi_1^*\alpha\big).$$
Recall that we identify $\Gamma$ with the current it represents.
Observe that if $\alpha$ is positive then $f^*(\alpha)$ and $f_*(\alpha)$
are positive closed $(p,p)$-currents which are smooth on a dense
Zariski open set and have no mass on analytic subsets of $X$. 
They are represented by forms with coefficients in $\Lc^1$.
For example, $f^*(\alpha)$ is smooth in $X\setminus P_1(f)$ and has no
mass on $P_1(f)$.
Moreover, if the positive closed $(p,p)$-forms $\alpha$ and
$\alpha'$ are cohomologous then
$\|f^*\alpha\|=\|f^*\alpha'\|$ and $\|f_*\alpha\|=\|f_*\alpha'\|$.
Define
$$\lambda_p(f):=\|f^*(\omega^p)\|=\int_X
f^*(\omega^p)\wedge\omega^{k-p}=\langle
\Gamma,\pi_2^*\omega^p\wedge\pi_1^*\omega^{k-p}\rangle = \int_X f_*(\omega^{k-p})\wedge\omega^p.$$
This integral can be computed cohomologically. It measures the norm of the linear operator $f^*$ acting
on the cohomology group $H^{p,p}(X,\C)$.

The following proposition shows that the sequence $c \lambda_p(f^n)$ is
sub-multiplicative, see also \cite{DinhSibony5}. 
Hence, $\lambda_p(f^n)^{1/n}$ converge to a constant $d_p(f)$. We call
$d_p(f)$ {\it the dynamical
degree of order $p$} of $f$. It is easy to check that $d_0(f)$ and $d_k(f)$ are the topological degrees
(i.e. the number of points in a generic fiber counted with multiplicities)
 of $\pi_{1|\Gamma}$ and $\pi_{2|\Gamma}$ and that $d_p(f^n)=d_p(f)^n$.

\begin{proposition} \label{prop_pullback_composition}
Let $f$ and $f'$ be two correspondences on $(X,\omega)$. Then, 
there exists a constant $c>0$ independent of $f$ and $f'$ such that 
$$\lambda_p(f\circ f')\leq c \lambda_p(f)\lambda_p(f').$$
\end{proposition}

We will need the following lemma.

\begin{lemma} \label{lemma_composition}
Let $\Omega\subset X\setminus P_1(f)$ and $\Omega'\subset X\setminus
P_1(f')$ be dense Zariski open subsets of $X$ such that
$f(\Omega)\subset X\setminus P_2(f)$ and
$f'(\Omega')\subset \Omega\setminus P_2(f')$. 
If $S$ is an arbitrary current on $X$, then
$$(f\circ f')^*_{|\Omega'}(S)={f'}^*_{|\Omega'}\big(f^*_{|\Omega}(S)\big).$$   
\end{lemma}
\proof
Let $U$ be a small neighbourhood of
a point in $\Omega'$. Since $\Omega'\cap P_1(f')=\varnothing$, the
restriction of $f'=\pi_2\circ (\pi_{1|\Gamma'})^{-1}$ to $U$ is
given by a family of biholomorphic maps $u_r:U\rightarrow
U_r\subset \Omega$. If $U$ is small enough, $f=\pi_2\circ
(\pi_{1|\Gamma})^{-1}$ restricted to each $U_r$ is given by a family
of biholomorphic maps $u_{rs}:U_r\rightarrow U_{rs}\subset X$. 
Hence $f\circ f'$ restricted to $U$ is given by the family
of biholomorphic maps $u_{rs}\circ u_r:U\rightarrow U_{rs}$. We
have
\begin{eqnarray*}
{f'}^*_{|U}\big(f^*_{|\Omega}(S)\big) & = & \sum_r u_r^*\big(\sum_s u_{rs}^*(S)\big)= \sum_{r,s}
u_r^*\big(u_{rs}^*(S)\big) \\
& = & \sum_{r,s}
(u_{rs}\circ u_r)^*(S) = (f\circ f')^*_{|U}(S).
\end{eqnarray*}
This implies the lemma. 
\endproof

\noindent
{\it Proof of Proposition \ref{prop_pullback_composition}.}
Observe that $(f\circ f')^*(\omega^p)$ is a positive closed current on
$X$. Moreover, $(f\circ f')^*(\omega^p)$ and $f'^*(f^*(\omega^p))$
are well defined and smooth outside an analytic set.
We obtain from Lemma \ref{lemma_composition}  
that these forms are equal 
on some Zariski open set $\Omega'$.
By Theorem \ref{th_regularization},
there exist positive closed  smooth $(p,p)$-forms $T_n$, converging to
a current $T\geq f^*(\omega^p)$, such that
$\|T_n\|\leq c\|f^*(\omega^p)\|=
c\lambda_p(f)$. Hence, there is another constant
$c>0$ such that 
$c\lambda_p(f)\omega^p-T_n$ is 
cohomologous to a smooth positive closed form for every $n$.
We have
$$\|{f'}^*_{|\Omega'} \big(f^*(\omega^p)\big)\|\leq\limsup_{n\rightarrow\infty} 
\|{f'}^*(T_n)\|\leq c\lambda_p(f)\|{f'}^*(\omega^p)\|=c\lambda_p(f)\lambda_p(f'). $$
Hence, since $(f\circ f')^*(\omega^p)$ has no mass on analytic sets,
$$\|(f\circ f')^*(\omega^p)\|=\|{f'}^*_{|\Omega'}\big(f^*(\omega^p)\big)\|\leq c\lambda_p(f)\lambda_p(f').$$
\endproof

\begin{remark} \rm
Let $A_{p,q}(f)$ denote the norm of $f^*$ on $H^{p,q}(X,\C)$. One can
prove as in \cite{Dinh1} that 
$$A_{p,q}(f)\leq c \sqrt{\lambda_p(f)\lambda_q(f)}$$
where $c>0$ is a constant independent of $f$. This inequality
and the Lefschetz fixed points formula allow to get an asymptotic
estimate of the number of periodic points of order $n$ of $f$ when
they are isolated. For example, if $d_k(f)$ is strictly larger than
the other dynamical degrees, this number is equal to $d_k(f)^n\big(1+o(1)\big)$. 
\end{remark}


\section{Entropy} \label{section_entropy}

We now define the topological entropy of $f$.
We call {\it $n$-orbit} of $f$ any sequence
$$(x_0,i_1,x_1,i_2,x_2,\ldots, x_{n-1},i_n, x_n)$$
where $x_0$, $\ldots$, $x_n$ are points of $X$ with $x_i\not\in I_1(f)$, 
and $i_1$, $\ldots$, $i_n$ are indices such that
$(x_{r-1},x_r)\in \Gamma_{i_r}$ for every $r$. 
Let $\Fc$ be a finite family of $n$-orbits of $f$.
We say that $\Fc$ is {\it $\epsilon$-separated} if for all distinct elements 
$$(x_0,i_1,x_1,i_2,x_2,\ldots, x_{n-1},i_n, x_n) \quad
\mbox{and} \quad (x'_0,i'_1,x'_1,i'_2,x'_2,\ldots, x'_{n-1},i'_n, x'_n)$$
of $\Fc$, we have either $i_r\not=i_r'$ or 
$\dist(x_r,x_r')>\epsilon$ for some $r$.
As we already explained, 
the indices $i_r$ allow  to count the multiplicities. When $\Gamma$ is irreducible, we always
have $i_r=i_r'$, then the indices $i_r$ in the definition of  $n$-orbit can be dropped. 
But since we are going to consider the graph of $f^n$, we cannot deal
only with the irreducible case. 

\begin{definition}[see \cite{Bowen, Gromov, Friedland}] \rm \label{def_entropy}
Define {\it the topological entropy} of $f$ by
$$h(f):=\sup_{\epsilon>0} \lim_{n\rightarrow\infty} \frac{1}{n}\log\max\big\{\# \Fc,\quad \Fc \mbox{ an }
\epsilon\mbox{-separated family of } n\mbox{-orbits of } f
\big\}.$$
\end{definition}

We say that the $n$-orbit
$\Oc=(x_0,i_1,x_1,i_2,x_2,\ldots, x_{n-1},i_n, x_n)$
is {\it regular} if for every $1\leq s\leq n$, $\Gamma_{i_s}$ is,
in a neighbourhood of $(x_{s-1}, x_s)$, a graph over each factor of $X^2$. Since
any $n$-orbit can be approximated by regular $n$-orbits, in
Definition \ref{def_entropy}, one can consider only regular orbits.
This is why we will consider only the extension by zero of all the currents
defined on a Zariski open set.

As observed in \cite{Gromov, Friedland}, $f$ is conjugated to a shift
$\sigma$ on the space $X_\infty$, the closure of the set of the infinite
orbits  $(x_0,i_1,x_1,i_2,\ldots, i_n,x_n,\ldots)\in X^\N$. It follows that
$h(f)=h(\sigma)$, and since $\sigma$ is continuous, one gets that $h(f^n)=nh(f)$ and
$h(f^{-1})= h(f)$.

Let $M=\{m_1,\ldots,m_s\}$, with $0\leq
m_1\leq m_2\leq \cdots \leq m_s$, be a multi-index. We define the graph $\Gamma_M$ of
$(f^{m_1},\ldots,f^{m_s})$ in $X^{s+1}$ as the closure of the set of
points $(x_0,x_{m_1},\ldots,x_{m_s})\in X^{s+1}$ associated to a 
regular $m_s$-orbit
$$\Oc=(x_0,i_1,x_1,i_2,x_2,\ldots, x_{m_s-1},i_{m_s}, x_{m_s}).$$
This is a holomorphic $k$-chain in $X^{s+1}$ 
where the multiplicity of a generic point $(x_0,x_{m_1},\ldots,x_{m_s})$ 
in $\Gamma_M$ is the number of the associated regular $m_s$-orbits
$\Oc$. If $M=\{n\}$ we obtain the graph $\Gamma_n$ of $f^n$ in $X^2$.
If $M=\{1,\ldots,n\}$, we obtain the graph $\Gamma_{[n]}$ of
$(f,f^2,\ldots, f^n)$ in $X^{n+1}$. Recall that
$$\lov(f)=\limsup_{n\rightarrow \infty} \log
\left(\volume(\Gamma_{[n]})^{1/n}\right)=
\limsup_{n\rightarrow \infty} \log \|\Gamma_{[n]}\|^{1/n}.$$
We divide the proof of Theorem
\ref{th_entropy} in two parts.

\medskip
\noindent
{\it Proof of the inequality.} 
We follow an idea due to Gromov \cite{Friedland}, see also \cite{Gromov}.
Let $\Fc$ be an $\epsilon$-separated family of regular 
$n$-orbits of $f$. We have to compare 
$\#\Fc$ with $\|\Gamma_{[n]}\|$. 
We associate to each element $\Oc=(x_0,i_1,\ldots,i_n,x_n)$ of $\Fc$ an open set 
$\Bc_\Oc\subset\Gamma_{[n]}$ which is the set of
the points
$$(x_0',\ldots,x_n') \in X^{n+1} \ \mbox{ with }(x'_{r-1},x'_r)\in \Gamma_{i_r} 
 \mbox{ and } \dist(x'_r,x_r)<\epsilon/2 \mbox{ for every } r.$$
Here, the distance between two points in $X^{n+1}$ is the maximum of
the distances between their projections on factors of $X^{n+1}$. 
Since $\Fc$ is $\epsilon$-separated, the balls $\Bc_\Oc$ are disjoint
(two balls with indices $(i_1,\ldots,i_n)\not=(i_1',\ldots,i_n')$ are considered as disjoint balls).
Hence, the total mass of all the $\Bc_\Oc$ is smaller
than $\|\Gamma_{[n]}\|$.

On the other hand, $\Bc_\Oc$ contains an analytic subset of dimension $k$ of the ball of diameter $\epsilon$ and 
of center $(x_0,\ldots, x_n)$ in $X^{n+1}$. A theorem of Lelong \cite{Lelong} 
implies that $\|\Bc_\Oc\|\geq c\epsilon^{2k}$ where
$c>0$ is a constant independent of $\epsilon$ and of $n$.
 
Let $\Pi_i:X^{n+1}\rightarrow X$ denote the canonical projections on
the factor of index $i$, $0\leq i\leq n$. Define
$\omega_i:=\Pi_i^*(\omega)$. 
We use for $X^{n+1}$ the  {\it canonical} K{\"a}hler form $\omega_0+\cdots+\omega_n$.
Then, the number of $\Bc_\Oc$, which is equal to $\#\Fc$, satisfies
$$\#\Fc\leq c^{-1}\epsilon^{-2k}\|\Gamma_{[n]}\|.$$
The inequality $h(f)\leq \lov(f)$
in Theorem \ref{th_entropy} follows from Definition \ref{def_entropy}.
\hfill $\square$

\medskip
\noindent
{\it Proof of the equality.} Recall that $\pi_1$,
$\pi_2:X^2\rightarrow X$ denote the canonical projections.
We have
$$\|\Gamma_n\|  =  \langle
\Gamma_n,(\pi_1^*\omega + \pi_2^*\omega)^k\rangle 
=\sum_{p=0}^k\left({k\atop p}\right) \langle \Gamma_n,
\pi_1^*\omega^{k-p}\wedge \pi_2^*\omega^p\rangle 
 =  \sum_{p=0}^k \left({k\atop p}\right)  \lambda_p(f^n).$$
Hence
\begin{eqnarray}
\max_{0\leq p\leq k} \lambda_p(f^n) \leq \|\Gamma_n\| \leq 2^k
\max_{0\leq p\leq k} \lambda_p(f^n). \label{eq_volume}
\end{eqnarray}
On the other hand, the projection of $\Gamma_{[n]}$ on the product
$X^2$ of the first and the last factors of $X^{n+1}$, 
is equal to $\Gamma_n$. 
It follows that $\|\Gamma_{[n]}\|\geq \|\Gamma_n\|$, hence $\lov(f)\geq \max \log d_p(f)$. 

For the other inequality, it is enough to show that
$\|\Gamma_{[n]}\|\lesssim n^k(\delta+\epsilon)^n$, where $\delta:=\max_p
d_p(f)$ and $\epsilon$ is a fixed constant.
We have
$$\|\Gamma_{[n]}\|=\langle \Gamma_{[n]},
(\omega_0+\cdots+\omega_n)^k\rangle
= \sum_{0\leq n_s\leq n} \langle \Gamma_{[n]}, \omega_{n_1}\wedge\ldots\wedge \omega_{n_k}\rangle.$$
We only need to prove that $\langle \Gamma_{[n]},
\omega_{n_1}\wedge\ldots\wedge \omega_{n_k}\rangle\leq
c(\delta+\epsilon)^n$, $c>0$. The following proposition  will be useful for that purpose.

\begin{proposition} \label{prop_volume}
There is a constant $c_s>0$ independent of the multi-index
$M=\{m_1,\ldots,m_s\}$, $0\leq m_1\leq\cdots\leq m_s$, such that 
$$\|\Gamma_M\|\leq c_s(\delta+\epsilon)^{m_s}.$$
\end{proposition}
\proof
The proof uses an induction on $s$. For $s=1$ we have
$\Gamma_M=\Gamma_{m_1}$, and the desired estimate follows from the
relation (\ref{eq_volume}).

Assume the proposition for $|M|=s-1$. We will prove it for $|M|=s$.
Let $\tau_1:X^{s+1}\rightarrow X^2$ be the canonical projection on the two
first factors and let $\tau_2:X^{s+1}\rightarrow X^{s}$ be the projection
on the $s$ last factors. Define $M':=\{m_2-m_1,\ldots,m_s-m_1\}$. We
will prove that
$\Gamma_M=\tau_1^{-1}(\Gamma_{m_1})\cap\tau_2^{-1}(\Gamma_{M'})$ in a
Zariski open set, then we will apply Corollary \ref{cor_intersection}.

Let $\Omega\subset X$ be the Zariski open set of all the
points $x_0\in X$ which admit $d_0(f)^{m_s}$ distinct regular
$m_s$-orbits, i.e. the maximal number of regular $m_s$-orbits.
Let $\Omega_{s+1}$ denote the Zariski open subsets of points in
$X^{s+1}$ whose projections on the first factor $X$ belong to
$\Omega$. Observe that $\Gamma_M\cap \Omega_{s+1}$ is Zariski dense in
$\Gamma_M$.
Hence, we only need to estimate $\|\Gamma_M\cap\Omega_{s+1}\|$. 

Consider a regular $m_s$-orbit
$\Oc:=(x_0,i_1,\ldots,i_{m_s},x_{m_s})$, $x_0\in\Omega$,  
associated to a point $z$ in $\Gamma_M \cap \Omega_{s+1}$. The point
$\tau_1(z)$ is associated to the regular $m_1$-orbit
$\Oc_1:=(x_0,i_1,\ldots,i_{m_1},x_{m_1})$, i.e. to a point in $\Gamma_{m_1}$.
The point
$\tau_2(z)$ is associated to the regular $(m_s-m_1)$-orbit
$\Oc_2:=(x_{m_1},i_{m_1},\ldots,i_{m_s},x_{m_s})$, i.e. to a point in
$\Gamma_{M'}$.
It follows that in $\Omega_{s+1}$, $\Gamma_M$ is the intersection of
$\tau_1^{-1}(\Gamma_{m_1})$ with $\tau_2^{-1}(\Gamma_{M'})$.

Let $\Omega_2$ denote the Zariski open subset of points in
$X^2$ whose projections on the first factor $X$ belong to
$\Omega$. The choice of $\Omega$ implies that in $\Omega_2$,
$\Gamma_{m_1}$ is locally a graph over the second factor $X$ of
$X^2$. It follows that $\tau_2$ restricted to
$\tau_1^{-1}(\Gamma_{m_1})\cap \Omega_{s+1}$  is locally biholomorphic.
Then, we can apply Corollary \ref{cor_intersection} and Remark
\ref{rk_intersection} to $\pi=\tau_2$,
and to  components of $\tau_1^{-1}(\Gamma_{m_1})$ and of
$\Gamma_{M'}$. We obtain
$$\|\Gamma_M\|=\|\Gamma_M\cap\Omega_{s+1}\|  \leq 
c\|\tau_1^{-1}(\Gamma_{m_1})\|\ \|\Gamma_{M'}\|
\leq  c'\|\Gamma_{m_1}\|\ \|\Gamma_{M'}\| $$
where $c$, $c'$ depend only on $(X,\omega)$ and on $s$.
The case $|M|=1$ and the case $|M|=s-1$ imply the result. 
\endproof

\noindent
{\it End of the proof of Theorem \ref{th_entropy}.} We will prove that  $\langle \Gamma_{[n]},
\omega_{n_1}\wedge\ldots\wedge \omega_{n_k}\rangle\leq
c(\delta +\epsilon)^n$, $c>0$, for $0\leq n_1\leq\cdots\leq n_k\leq n$.
Let $\Pi:X^{n+1}\rightarrow X^{k+1}$ be the canonical projection on the
product of factors with indices $0$, $n_1$, $\ldots$, $n_k$. 
We show that $\Pi$ defines a map of topological degree 
$d_0(f)^{n-n_k}$ between $\Gamma_{[n]}$ and  $\Gamma_M$, where
$M:=\{n_1,\ldots,n_k\}$.

Observe that if we fix a generic orbit
$\Oc':=(x_0,i_1,\ldots,i_{n_k},x_{n_k})$ there are 
$d_0(f)^{n-n_k}$ choices for $\Oc'':=(i_{n_k+1},x_{n_k+1},\ldots,i_n,x_n)$ such
that 
$\Oc:=(\Oc',\Oc'')$ is a point in $\Gamma_{[n]}$. By
definition, $\Oc'$
corresponds to a point in $\Gamma_M$. Hence,  $\Pi$ defines a map of topological degree 
$d_0(f)^{n-n_k}$ between $\Gamma_{[n]}$ and  $\Gamma_M$.

If $\widetilde\omega$ denotes the canonical K{\"a}hler form on
$X^{k+1}$, then
$$\langle \Gamma_{[n]}, \omega_{n_1}\wedge\ldots\wedge
\omega_{n_k}\rangle  \leq   \langle
\Gamma_{[n]},\Pi^*(\widetilde\omega^k)\rangle 
 =  d_0(f)^{n-n_k} \langle \Gamma_M,\widetilde\omega^k\rangle 
  =  d_0(f)^{n-n_k} \|\Gamma_M\|.$$
By Proposition \ref{prop_volume},  $\|\Gamma_M\|\leq c_k
(\delta+\epsilon)^{n_k}$. The desired estimate follows.
\hfill $\square$


\section{Entropy on a subvariety} \label{partial_entropy}

Let $Y\subset X$ be an analytic subset of pure dimension $m$ or more
generally a holomorphic $m$-chain. Assume
that for a generic point $x_0\in Y$ the sets $f^n(x_0)$ do not
intersect $I_1(f)$ for any $n\geq 0$. Such a point admits $n$-orbits. 
We define the entropy $h(f,Y)$ of $f$ on $Y$ as in
Definition \ref{def_entropy} but we only consider  the orbits 
$\Oc=(x_0,i_1,\ldots,i_n,x_n)$
starting from a point $x_0\in Y$. 
Define the holomorphic chain $\Gamma^Y_{[n]}$ as the closure in
$X^{n+1}$ of the set of $n$-orbits $\Oc=(x_0,i_1,x_1,\ldots,i_n,x_n)$ with
$x_0\in Y$ generic, and
$$\lov(f,Y):=\limsup_{n\rightarrow \infty} \log\big(\volume
(\Gamma^Y_{[n]} )^{1/n}\big).$$

We have the following result which generalizes Theorem \ref{th_entropy}.

\begin{theorem} \label{th_entropy_bis}
Let $Y$ be as above. Assume that all the orbits starting from a generic
point $x_0\in Y$ are regular \footnote{this hypothesis is in fact not
  necessary, but the proof for the general case needs a
  theory of intersection of currents that we will develop in a future
  work; for meromorphic maps this hypothesis is clearly satisfied}.
Then
$$h(f,Y)\leq \lov(f,Y)\leq \max_{0\leq p \leq m} \log d_p(f).$$
\end{theorem}

Such an estimate should be useful in the study of dimensional entropies
and Lyapounov exponents. We refer to Newhouse \cite{Newhouse} and Buzzi
\cite{Buzzi} for this purpose. 
The proof 
uses the same idea as in Theorem \ref{th_entropy}. The first inequality is left to the reader.
For the second inequality, in order to estimate $\volume(\Gamma_{[n]}^Y)$,
it is sufficient to apply Proposition
\ref{prop_volume_bis} below for $T=Y$. Proposition
\ref{prop_volume_bis} is more general than Proposition
\ref{prop_volume}. However, we keep Proposition
\ref{prop_volume} because its proof contains a useful geometric argument.

Let $M=\{m_1,\ldots,m_s\}$, $0\leq m_1\leq\cdots\leq m_s$, be a
multi-index. Let $\widetilde \Gamma_M$ denote the largest Zariski open subset of 
$\Gamma_M$ which is locally a graph over the first factor of $X^{s+1}$. 
Define $u:\widetilde\Gamma_M\rightarrow X$ the canonical projection on the first
factor and $\delta_m:=\max_{0\leq p\leq m}d_p(f)$.

\begin{proposition} \label{prop_volume_bis}
There is a constant $c_s>0$ independent of $M$ such that if $T$ is a
positive closed $(k-m,k-m)$-current on $X$ then $u^*(T)$ defines
a positive closed current of bidimension $(m,m)$ on $X^{s+1}$ with
$$\|u^*(T)\|\leq c_s\|T\|(\delta_m+\epsilon)^{m_s}.$$
\end{proposition}
\proof
By Skoda's theorem \cite{Skoda}, the trivial extension of $u^*(T)$
is positive and closed in $X^{s+1}$ provided that $u^*(T)$ has
finite mass. So, it is enough to estimate 
$\|u^*(T)\|$. By Theorem \ref{th_regularization}, the case where $T$ is smooth implies
the general case. Hence, we can assume $T$ smooth.
The proof uses an induction on $s$. 

For $s=1$ we have
$\Gamma_M=\Gamma_{m_1}$. We need to show that $\langle u^*(T),
\omega_0^r\wedge\omega_1^{m-r}\rangle \lesssim
\|T\|(\delta_m+\epsilon)^{m_1}$ for $0\leq r\leq m$. 
Choose a constant $c>0$, independent of $T$, such that
$c\|T\|\omega^{k-m+r}-T\wedge\omega^r$ is cohomologous to a smooth positive
closed form. Hence
\begin{eqnarray*}
\langle u^*(T)\wedge \omega_0^r,\omega_1^{m-r}\rangle & = &  \langle
u^*(T\wedge\omega^r),\omega_1^{m-r}\rangle \\
& \leq &  c\|T\| \langle
\omega_0^{k-m+r}\wedge \Gamma_{m_1}, \omega_1^{m-r}\rangle \\
& = &  c\|T\|\lambda_{m-r}(f^{m_1}).
\end{eqnarray*}
The desired estimate follows.

Now, assume the inequality for $|M|=s-1$ and for arbitrary $T$, smooth or
not. We will prove it for $|M|=s$
and for $T$ smooth.
We have to estimate $\langle
u^*(T),\omega_0^{r_0}\wedge \ldots \wedge \omega_s^{r_s}\rangle$
with $r_0+\cdots+r_s=m$. Since this integral is equal to $\langle
u^*(T\wedge \omega^{r_0}),\omega_1^{r_1}\wedge \ldots\wedge
\omega_s^{r_s}\rangle$, we can replace $T$ by $T\wedge \omega^{r_0}$
and assume that $r_0=0$. 
Now, let $\tau_2$ be, as in Proposition \ref{prop_volume}, the
projection from $X^{s+1}$ on the last $s$ factors. Define 
 $\Theta:=\widetilde\omega^m$ where $\widetilde\omega$ is
 the canonical K{\"a}hler form on $X^s$. It is enough to estimate 
$\langle u^*(T),\tau_2^*(\Theta)\rangle = \langle \tau_{2*}(u^*(T)),
\Theta \rangle$.

Observe that $\tau_{2*}(u^*(T))$ is supported in
$\Gamma_{M'}$ (see Proposition \ref{prop_volume}) 
and has no mass on analytic subsets of $\Gamma_{M'}$,
since $T$ is smooth. 
Let $\widetilde\Gamma_{M'}$ denote the largest Zariski open subset of 
$\Gamma_{M'}$ which is locally a graph over the first factor $X$ of
$X^{s}$
and let $u':\widetilde\Gamma_{M'}\rightarrow X$ denote the canonical
projection on this factor.
We will prove as in Lemma
\ref{lemma_composition} that $\tau_{2*}(u^*(T))={u'}^*(f^{m_1}_*(T))$
on a Zariski open set of
$\Gamma_{M'}$. We first assume this  and complete the proof.

The case $s=1$ implies that 
$$\|f^{m_1}_*(T)\| =\big\|\pi_{2*}\big(\pi_1^*(T)\wedge\Gamma_{m_1}\big)\big\|
\leq \|\pi_1^*(T)\wedge\Gamma_{m_1}\|\lesssim
\|T\|(\delta_m+\epsilon)^{m_1}.$$ 
The case $|M|=s-1$, applied to $M'$ and to $f^{m_1}_*(T)$, yields 
$$\big\langle \tau_{2*}(u^*(T)),\Theta\big\rangle = 
\big\langle {u'}^*(f^{m_1}_*(T)),\Theta \big\rangle \lesssim \|f^{m_1}_*(T)\|
(\delta_m+\epsilon)^{m_s-m_1}\lesssim \|T\| (\delta_m+\epsilon)^{m_s}.$$
It follows that $\langle u^*(T),\tau_2^*(\Theta)\rangle \lesssim \|T\|
(\delta_m+\epsilon)^{m_s}$ which implies the result.

Now, we prove the identity $\tau_{2*}(u^*(T))={u'}^*(f^{m_1}_*(T))$ on
a Zariski open set of $\Gamma_{M'}$. Let $U$ be
a small neighbourhood of a generic point in $\Gamma_{M'}$. Then $u'$ defines a
biholomorphic map between $U$ and an open set $V\subset X$. If $U$ is
small enough, $f^{-m_1}$ restricted to $V$ is given by a family of
biholomorphic maps $u_r:V\rightarrow V_r\subset X$. 
Observe that a generic point $(x,z)\in \Gamma_M$ is sent by $\tau_2$ to $z\in\Gamma_{M'}$
if and only if $x$ is sent by $f^{m_1}$ to $u'(z)$.
Then, 
$\tau_{2|U}^{-1}$ is given by a family of biholomorphic maps between $U$
and the open sets 
$$U_r:=\big\{\big(u_r(u'(z)),z\big)\in X\times X^s,\ z\in U\big\}.$$
These maps, by definition of $\tau_2$, are equal to 
$z\mapsto \big(u_r(u'(z)),z\big)$.
From the definition of $u$, we deduce that $u \circ \tau_{2|U}^{-1}$
is given by the family of the
biholomorphic maps $u_r\circ u':U \rightarrow V_r$.
Hence, on $U$, we have
$$\tau_{2*}(u^*(T))=\sum_r (u_r\circ {u'})^* (T)={u'}^* \Big(\sum_r
u_r^*(T) \Big) = {u'}^*(f^{m_1}_*(T)).$$
This implies the result.
\endproof


\section{Julia and Fatou sets} \label{julia_fatou}

We discuss here  a notion of Julia and Fatou sets for
correspondences. Let $B_x(r)$ denote the ball of center $x$ and of
radius $r$. The following function, which describes the local growth
of volume of graph,  has strong links with the Julia and
Fatou sets :
$$\Phi(x):=\inf_{r>0} \limsup_{n\rightarrow\infty} {1\over
  n}\log\volume \big(\Gamma_n\cap\pi_1^{-1}(B_x(r))\big)$$
(we can also consider $\Gamma_{[n]}$ instead of $\Gamma_n$). 
Since $\pi_1$ restricted to $\Gamma_n$ has topological degree
$d_0(f)^n$, we have $\Phi(x)\geq \log d_0(f)$. Proposition \ref{prop_volume}
implies that $\Phi(x)\leq \max_p \log d_p(f)$. 
It is easy to check that the function $\Phi$ is upper semi-continuous. 
We can study the sets
$\{\Phi< \delta\}$ and $\{\Phi\geq\delta\}$ as analogues of Fatou and
Julia sets. It is likely that ergodic invariant measures of
maximal entropy, if they exist, are supported on the set where $\Phi$
take the maximal value. Consider some examples.

\begin{example} \rm
Let $f:\P^k\rightarrow\P^k$ be a holomorphic map of algebraic degree
$d\geq 2$. It is well known that $d^{-pn} f^{n*}(\omega^p_\FS)$ converges
to $T^p$. Here, 
$\omega_\FS$ denotes the Fubini-Study form on $\P^k$ and  $T$ denotes 
the Green $(1,1)$-current of $f$.
 The volume of
$\Gamma_n\cap\pi_1^{-1}(B_x(r))$ is the sum over $p$ of the integrals 
$\langle f^{n*}(\omega_\FS^p)_{|B_x(r)},\omega_\FS^{k-p}\rangle$ on
$B_x(r)$. One estimates these integrals using the speed of convergence
of $d^{-n}f^{n*}(\omega_\FS)$ toward $T$, see \cite{Sibony, deThelin,
  Dinh2}, and one deduces that
$\Phi(x)=\log d^p$ if $x\in\supp(T^p)\setminus\supp(T^{p+1})$. The support
$\supp(T^p)$ of $T^p$ and its complement are the Julia and the Fatou  sets of order $p$ 
associated to $f$. For $p=1$, one
obtain the classical Fatou and Julia sets, see
\cite{FornaessSibony}. The function $\Phi$ takes only $k+1$ values and
the set $\{\Phi\geq \log d^p\}$ supports the invariant  positive closed 
current $T^p$.
\end{example}

The following trivial example shows that, in general, Fatou and Julia sets cannot
be characterized only by the values of $\Phi$.

\begin{example} \rm
Consider $f:\P^1\rightarrow \P^1$ given by $z\mapsto 2z$ where $z$ is
an affine coordinate. Then $\Phi(x)=0$ everywhere but the family
$(f^n)_{n\geq 0}$ is locally equicontinuous except at 0.
The limit of $\Gamma_n$ contains a singular fiber $\pi_1^{-1}(0)$ as
component. Taking a product of $f$ with other holomorphic maps gives
analogous examples in any dimension with positive entropy.
\end{example}

One sees in the example below that the meromorphic case is quite more
delicate.

\begin{example} \rm
Let $f:\P^2\rightarrow\P^2$ be the meromorphic map given by $(z,w)\mapsto
(z^{-d}, w^{-d})$, $d\geq 2$, where $(z,w)$ denotes affine coordinates
of $\P^2$. Using the fact that $f^2(z,w)=(z^{d^2},w^{d^2})$, we
obtain that $\Phi(0)=0$; but $0$ is a point of indeterminacy of $f$.
\end{example}

Now define
$$\Psi(x):=\limsup_{r\rightarrow 0} \limsup_{n\rightarrow\infty}
{\volume\big(\Gamma_n \cap \pi_1^{-1}(B_x(r))\big) \over r^{2k} d_0(f)^n}.$$
It is left to the reader to check that if $f$ is a holomorphic
endomorphism of $\P^k$
then
$\{\Psi<\infty\}$ and $\{\Psi=\infty\}$ are the Fatou and Julia sets of $f$.


Tien-Cuong Dinh \hfill Nessim Sibony\\
Institut de Math{\'e}matique de Jussieu \hfill  Math{\'e}matique - B{\^a}timent 425 \\
Plateau 7D, Analyse Complexe \hfill  UMR 8628\\
175 rue du Chevaleret \hfill  Universit{\'e} Paris-Sud \\
75013 Paris, France \hfill 91405 Orsay, France \\
{\tt dinh@math.jussieu.fr} \hfill {\tt nessim.sibony@math.u-psud.fr} \\

\end{document}